Wael Bahsoun†, Igor V. Evstigneev*† and Michael I. Taksar‡

†Economics Department, University of Manchester, Oxford Road, Manchester M13 9PL, UK
‡ Mathematics Department, University of Missouri, Columbia, MO 65211, USA

The paper examines random dynamical systems related to the classical von Neumann and Gale models of economic growth. Such systems are defined in terms of multivalued operators in spaces of random vectors, possessing certain properties of convexity and homogeneity. A central role in the theory of von Neumann-Gale dynamics is played by a special class of paths called rapid (they maximize properly defined growth rates). Up to now the theory lacked quite satisfactory results on the existence of such paths. This work provides a general existence theorem holding under assumptions analogous to the standard deterministic ones. The result solves a problem that remained open for more than three decades.

*Keywords*: Random dynamical systems, convex multivalued operators, von Neumann–Gale model, rapid paths, Dynkin's and Radner's dualities

*AMS Subject Classification*: 37H99, 37H15, 91B62, 91B28

## 1. Introduction

Von Neumann-Gale dynamical systems are defined in terms of multivalued operators possessing certain properties of convexity and homogeneity. States of such systems are represented by elements of convex cones $X_t$ ($t = 0, 1, ...$) in linear spaces. Possible one-step transitions from one state to another are described in terms of given operators $A_t(x)$, assigning to each $x \in X_{t-1}$ a convex subset $A_t(x) \subseteq X_t$. It is assumed that the graphs $Z_t := \{(x, y) \in X_{t-1} \times X_t : y \in A_t(x)\}$ of the operators $A_t(x)$ are convex cones. Paths of the von Neumann-Gale dynamical system are sequences $x_0 \in X_0$, $x_1 \in X_1, ...$ such that $x_t \in A_t(x_{t-1})$. The classical, deterministic theory of such dynamics was originally aimed at the modeling of economic growth (von Neumann [26], Gale [14]). First attempts to build a stochastic generalization of this theory were undertaken in the 1970s by Dynkin [5–7], Radner [19; 20] and others. However, the initial attack on the problem left many questions unanswered. Substantial progress was made only in the 1990s—see the survey in [10].

It has recently been observed [2] that stochastic analogues of von Neumann-Gale systems provide a natural and convenient framework for financial modeling (asset pricing and hedging under transaction costs). This observation gave a new momentum to studies in the field and posed new interesting questions. It also revived interest in old unsolved problems.

A central role in the stochastic theory of von Neumann-Gale dynamics is played by a special class of paths called *rapid* (they maximize properly defined growth rates). In spite of the current progress achieved, the theory lacked quite satisfactory results on the existence of such paths. The results available up to now established their existence under certain conditions which also guaranteed their stability (turnpike property)—see [13]. These conditions seemed to be too restrictive in the context of economic models and did not cover a number of new examples arising in financial applications. In this paper, we fill this gap by establishing the existence of rapid paths under assumptions fully analogous to the standard deterministic ones. The result obtained is applicable to both old models related to economics and new ones, coming from finance. It gives a final solution to a problem that remained open for more than three decades.

*Corresponding author. Email: Igor.Evstigneev@manchester.ac.uk

This study is connected with the paper [11], where closely related questions are considered for stationary systems. Some of the techniques used here (convex duality) are similar to those in [11]. However, instead of the method of elimination of randomization (Dvoretzky, Wald and Wolfowitz [4]), we use in the existence proof a compactness principle based on a version of the multidimensional Fatou lemma (Schmeidler [24]).

In this paper we consider a class of stochastic von Neumann-Gale dynamical systems defined in terms of a probability space $(\Omega, \mathcal{F}, P)$, a filtration $\mathcal{F}_0 \subseteq \mathcal{F}_1 \subseteq ... \subseteq \mathcal{F}$ and a sequence of random closed convex cones $G_t(\omega) \subseteq \mathbb{R}^n_+ \times \mathbb{R}^n_+$ ($t = 1, 2, ...$). For each $t$, it is assumed that the cone $G_t(\omega)$ depends $\mathcal{F}_t$-measurably[1] on $\omega$, and for all $t$ and $\omega$, the following conditions hold:

(**G.1**) for any $a \in \mathbb{R}^n_+$, the set $\{b : (a, b) \in G_t(\omega)\}$ is non-empty;

(**G.2**) the set $G_t(\omega)$ is contained in $\{(a, b) : |b| \leq M_t|a|\}$, where $M_t$ is a constant independent of $\omega$ (we write $|\cdot|$ for the sum of the absolute values of the coordinates of a vector);

(**G.3**) there exist a strictly positive constant $\gamma_t > 0$ and a pair of essentially bounded vector functions $(\hat{a}_{t-1}(\omega), \hat{b}_t(\omega))$ such that $(\hat{a}_{t-1}(\omega), \hat{b}_t(\omega)) \in G_t(\omega)$ for all $\omega$ and $\hat{b}_t(\omega) \geq \gamma_t e$, where $e = (1, ..., 1)$ (all inequalities between vectors are understood coordinatewise).

Denote by $\mathcal{X}_t$ the set of $\mathcal{F}_t$-measurable essentially bounded vector functions $x(\omega)$ with values in $\mathbb{R}^n_+$ and put

$$Z_t := \{(x, y) \in \mathcal{X}_{t-1} \times \mathcal{X}_t : (x(\omega), y(\omega)) \in G_t(\omega) \text{ (a.s.)}\}$$

($t = 1, 2, ...$). We will consider the von Neumann-Gale dynamical system defined by the multivalued operators

$$A_t(x) := \{y \in \mathcal{X}_t : (x, y) \in Z_t\} \ (x \in \mathcal{X}_{t-1}, \ t = 1, 2, ...).$$

A finite or infinite sequence $(x_t(\omega))_{t=0}^N$, $x_t \in \mathcal{X}_t$, $1 \leq N \leq \infty$, is called a *path* (*trajectory*) in the system under consideration if $x_t \in A_t(x_{t-1})$ or equivalently if

$$(x_{t-1}(\omega), x_t(\omega)) \in G_t(\omega) \text{ (a.s.)}$$

for all $t \geq 1$ for which $x_t$ is defined.

Define

$$G_t^\times(\omega) := \{(c, d) \geq 0 : db - ca \leq 0 \text{ for all } (a, b) \in G_t(\omega)\}, \tag{1}$$

where $ca$ and $db$ denote the scalar products of the vectors. Let $\mathcal{P}_t$ denote the set of $\mathcal{F}_t$-measurable vector functions $p(\omega)$ with values in $\mathbb{R}^n_+$ such that $E|p(\omega)| < \infty$. A *dual path* (*dual trajectory*) is a finite or infinite sequence $p_1(\omega), p_2(\omega), ...$ such that $p_t \in \mathcal{P}_t$ ($t \geq 1$) and

$$(p_t(\omega), E_t p_{t+1}(\omega)) \in G_t^\times(\omega) \text{ (a.s.)} \tag{2}$$

for all $t \geq 1$ for which $p_{t+1}$ is defined. We write $E_t(\cdot) = E(\cdot|\mathcal{F}_t)$ for the conditional expectation given $\mathcal{F}_t$. A dual path $(p_1, ..., p_{N+1})$ is said to *support* a path $(x_0, ..., x_N)$ if

$$p_t x_{t-1} = 1 \text{ (a.s.)} \tag{3}$$

---

[1] A closed set $G(\omega)$ in a metric space is said to depend measurably on $\omega$ if the distance to this set from each point in the space is a measurable function of $\omega$.

for $t = 1, ..., N + 1$ (for infinite paths (3) should hold for all $t \geq 1$). A trajectory is called *rapid* if there exists a dual trajectory supporting it. The term "rapid" is motivated by the fact that

$$\frac{E_t(p_{t+1}y_t)}{p_t y_{t-1}} \leq \frac{E_t(p_{t+1}x_t)}{p_t x_{t-1}} = 1 \text{ (a.s.)}$$

for each path $y_0, y_1, ...$ with $p_t y_{t-1} > 0$ (see (2) and (3)). This means that the path $x_0, x_1, ...$ maximizes the conditional expectation of the *growth rate* at each time $t$, the maximum being equal to 1. Growth rates are measured in terms of the random linear functions $p_t a$. In the applications, states of the system $x_t$ constituting trajectories can represent, e.g., commodity vectors or portfolios of assets; then $p_t$ are interpreted as price vectors. For a comprehensive discussion of the applications of von Neumann-Gale systems in economics and finance see [10] and [2].

The main result of the paper is the following theorem.

**Theorem 1.** *Let $x_0(\omega)$ be a vector function in $\mathcal{X}_0$ such that $\delta e \leq x_0(\omega) \leq De$ for some constants $0 < \delta \leq D$. (i) For each $N \geq 1$, there exists a finite rapid path of length $N$ with initial state $x_0(\omega)$. (ii) There exists an infinite rapid path with initial state $x_0(\omega)$.*

Assertion (i) of this theorem will be proved in the next section. Assertion (ii) will be proved in Section 3.

## 2. Dynkin's and Radner's dualities

To prove assertion (i) of Theorem 1 we will use some results of previous work, dealing with a somewhat different definition of dual paths. In the definition given in Section 1, we follow essentially the approach of Dynkin [5; 6]. Another version of this concept was introduced by Radner (see, e.g., [21; 22]). Dynkin and Radner both dealt with somewhat different models, and when referring to their work, we adjust their considerations for the present context. According to Radner's approach (as applied in our context), a dual path is defined as a sequence of non-negative measurable vector functions $q_0(\omega), q_1(\omega), ...$ such that $q_t \in \mathcal{P}_t$ and

$$E_{t-1}(q_t v) \leq q_{t-1} u \text{ (a.s.)} \tag{4}$$

for all $(u, v) \in Z_t$. We will refer to such sequences as *R-dual paths*. Since $G_t(\omega)$ is a cone, the set $Z_t$ contains with each pair of vector functions $(u, v)$ all pairs $(\lambda u, \lambda v)$, where $\lambda(\omega) \geq 0$ is any bounded $\mathcal{F}_{t-1}$-measurable scalar function. This implies that condition (4) is equivalent to

$$E(q_t v) \leq E(q_{t-1} u), \ (u, v) \in Z_t. \tag{5}$$

An $R$-dual path $q_0, q_1, ...$ is said to *support* a path $x_0, x_1, ...$ if

$$q_t(\omega) x_t(\omega) = 1 \text{ (a.s.)}. \tag{6}$$

Radner's approach appears to be more natural in the economic applications: $R$-dual paths admit a clear interpretation as price systems supporting optimal trajectories of economic dynamics. In this connection, most of the studies on stochastic von Neumann-Gale dynamical systems have dealt with what we call here $R$-dual paths. However, it has been recently shown [2] that Dynkin's approach is more suitable for the applications related to finance. In the present work, this approach serves as a powerful mathematical tool, playing a substantial role in the proof of the main existence theorem. Dealing with the notion of a dual path defined in accordance with Dynkin's approach, we employ a number of known results about $R$-dual paths and the following fact.

**Theorem 2.** *A path $(x_t)_{t=0}^N$ ($1 \leq N \leq \infty$) is supported by an R-dual path if and only if it is supported by a dual path.*

In the proof of this theorem we follow, with some necessary changes, the line of arguments in [11], where an analogous result was obtained in a stationary setting.

*Proof of Theorem 2.* "If". Let $x_0, x_1, ...$ be a (finite or infinite) path and $p_1, p_2, ...$ a dual path supporting it, i.e. satisfying (2) and (3). Define $q_t = E_t p_{t+1}$. By virtue of (2), $[E_t p_{t+1}] v \leq p_t u$ (a.s.) for all $(u, v)$ in $Z_t$. By taking the conditional expectation $E_{t-1}$ of both sides of the last inequality, we obtain (4). Condition (6) holds as a consequence of (3), and so $q_0, q_1, ...$ is an $R$-dual path supporting $x_0, x_1, ...$.

"Only if." To prove the converse assertion consider an $R$-dual path $q_0, q_1, ...$ supporting a path $x_0, x_1, ...$. We know that inequality (5) hold. By virtue of Proposition 1 we prove below, this implies the existence of an $\mathcal{F}_t$-measurable vector function $g_t$ such that $E|g_t| < \infty$, $E_{t-1} g_t = 0$ (a.s.) and

$$E[q_t(\omega)v(\omega) - q_{t-1}(\omega)u(\omega)] - E g_t(\omega) u(\omega) \leq 0 \tag{7}$$

for all $(u, v) \in \mathcal{X}_t \times \mathcal{X}_t$ such that $(u(\omega), v(\omega)) \in G_t(\omega)$ (a.s.). We will denote the class of such functions $(u(\omega), v(\omega))$ by $W_t$.

Define $p_t(\omega) = q_{t-1}(\omega) + g_t(\omega)$, $t \geq 1$. We claim that $p_1, p_2, ...$ is dual path supporting $x_0, x_1, ...$. To prove this, from (7) we obtain

$$E[q_t(\omega)v(\omega) - p_t(\omega)u(\omega)] \leq 0 \tag{8}$$

for all $(u, v) \in W_t$. By virtue of (**G.1**) and (**G.2**), for any bounded $\mathcal{F}_t$-measurable function $u(\omega) \geq 0$, there exists bounded $\mathcal{F}_t$-measurable $v(\omega) \geq 0$ such that $(u(\omega), v(\omega)) \in G_t(\omega)$ a.s. (we use here Aumann's measurable selection theorem). In view of (8), this implies $E p_t(\omega) u(\omega) \geq E q_t(\omega) v(\omega) \geq 0$ for any bounded $\mathcal{F}_t$-measurable $u(\omega) \geq 0$. Consequently, $p_t(\omega) \geq 0$ (a.s.), and so $p_t \in \mathcal{P}_t$.

Denote by $\hat{G}_t(\omega)$ the set of those $(a, b) \in G_t(\omega)$ which satisfy $|a| \leq 1$ and by $\hat{W}_t$ the class of all $\mathcal{F}_t$-measurable functions $(u(\omega), v(\omega))$ such that $(u(\omega), v(\omega)) \in \hat{G}_t(\omega)$ (a.s.). By virtue of (**G.2**) all such $(u, v)$ are essentially bounded, and so $\hat{W}_t \subseteq W_t$. Inequality (8), holding for all $(u, v) \in \hat{W}_t$, implies that for almost all $\omega$,

$$q_t(\omega) b - p_t(\omega) a \leq 0 \text{ for all } (a, b) \in G_t(\omega) \text{ with } |a| \leq 1 \tag{9}$$

(see Proposition A.3 in the Appendix). We can drop the constraint $|a| \leq 1$ in (9) because $G_t(\omega)$ is a cone. Since $E_t p_{t+1} = E_t(q_t + g_{t+1}) = q_t$ (a.s.), we have that with probability one

$$(E_t p_{t+1}) b - p_t a \leq 0 \text{ for all } (a, b) \in G_t(\omega),$$

which yields $(p_t, (E_t p_{t+1})) \in G_t^\times(\omega)$ (a.s.). Thus $p_1, p_2, ...$ is a dual path.

By setting $a = x_{t-1}(\omega)$ and $b = x_t(\omega)$ in (9), we get $1 = q_t(\omega) x_t(\omega) \leq p_t(\omega) x_{t-1}(\omega)$ a.s. for all $t \geq 1$ (see (6)). On the other hand,

$$E p_t x_{t-1} = E E_{t-1}(p_t x_{t-1}) = E[(E_{t-1} p_t) x_{t-1}] = E q_{t-1} x_{t-1} = 1, \ t \geq 1.$$

Therefore $p_t x_{t-1} = 1$ (a.s.), which proves (3). Thus $p_1, p_2, ...$ is a dual path supporting $x_0, x_1, ...$. □

In the class $W_t$ of vector functions $(u, v)$ defined above, consider the subclass $W^{t-1}$ consisting of those $(u, v)$ for which there exists an $\mathcal{F}_{t-1}$-measurable function $u'(\omega)$ coinciding with $u(\omega)$ (a.s.). Clearly, those and only those $(u, v)$ in $W_t$ belong to $W^{t-1}$ which satisfy $u = E_{t-1} u$ (a.s.). Put $r := q_t \in \mathcal{P}_t$ and $q := q_{t-1} \in \mathcal{P}_{t-1}$. From (5), we get

$$F(u, v) := E[r(\omega) v(\omega) - q(\omega) u(\omega)] \leq 0 \tag{10}$$

for all $(u, v) \in W^{t-1}$. In the course of the proof of Theorem 2, we used the following fact.

**Proposition 1.** *If inequality (10) holds for all $(u,v) \in W^{t-1}$, then there exists an $\mathcal{F}_t$-measurable vector function $g(\omega)$ such that $E|g| < \infty$, $E_{t-1}g(\omega) = 0$ (a.s.) and inequality*

$$F(u,v) - Eg(\omega)u(\omega) \leq 0$$

*holds for all $(u,v) \in W_t$.*

Inequality (10) says that the maximum of the functional $F(u,v)$ on the set of $(u,v) \in W_t$ satisfying the "non-anticipativity constraint" $u = E_{t-1}u$ (a.s.) is zero. The function $g(\omega)$ plays the role of a Lagrange multiplier relaxing this constraint. Lagrange multipliers of this kind were first considered by Rockafellar and Wets [23]. In the proof below, we use techniques outlined in [9].

*Proof of Proposition 1.* Denote by $L^k_\infty(t)$ the Banach space of essentially bounded measurable vector functions $x(\omega)$ with values in $\mathbb{R}^k$. We will regard functions $w = (u,v)$ in $W_t$ as elements of $L^{2n}_\infty(t)$. Consider the operator $B$ transforming $w = (u,v) \in W_t$ into $Bw := u - E_{t-1}u$. This is a continuous linear operator of $L^{2n}_\infty(t)$ into $L^n_\infty(t)$. We can characterize $W^{t-1}$ as the set of those $w \in W_t$ for which $Bw = 0$. According to (10), the function $\bar{w} := 0$ is a solution to the problem of maximization of the linear functional $F(w)$ on the convex set $W_t$ subject to the linear constraint $Bw = 0$. To analyze this problem we will use a version of the Kuhn-Tucker theorem formulated in the Appendix as Proposition A.4. We will apply this result to the Banach spaces $D_1 := L^{2n}_\infty(t)$ and $D_2 := L^n_\infty(t)$. Note that the image $H := B(D_1)$ is a closed subspace of $L^n_\infty(t)$ because it consists of those $y \in L^n_\infty(t)$ for which $E_{t-1}y = 0$ (a.s.). To verify condition (*) of Proposition A.4 consider the set $\tilde{W}_t$ consisting of those $(u,v) \in W_t$ which satisfy $0 \leq u \leq e$ (a.s.). Observe that $0 \in B(\tilde{W}_t)$ because $(0,0) \in \tilde{W}_t$ and the set $\{y \in H : ||y||_\infty < 1/2\}$ is contained in $B(\tilde{W}_t)$. Indeed, for any $y \in H$ with $||y||_\infty < 1/2$ define $u := e/2 + y$. Then $u \in \mathcal{X}_t$ and, by using (**G.1**), (**G.2**) and the measurable selection theorem, we can find a $v \in \mathcal{X}_t$ such that $(u(\omega), v(\omega)) \in G_t(\omega)$ (a.s.). We have $w : (u,v) \in \tilde{W}_t$ because $0 \leq u \leq e$ (a.s.), and $Bw = B(e/2 + y, v) = y - E_{t-1}y = y$ (a.s.) as $y \in H$, and so $E_{t-1}y = 0$ (a.s.). Thus $B(\tilde{W}_t)$ contains the open neighbourhood $\{y \in H : ||y||_\infty < 1/2\}$ of 0 in the space $H$. It remains to observe that $F(w)$ is bounded on $\tilde{W}_t$ because $\tilde{W}_t$ is a bounded set in $L^{2n}_\infty(t)$ and the functions $r$ and $q$ involved in the definition of $F$ (see (10)) are integrable.

By virtue of Proposition A.4, there exists a functional $\pi \in [L^{2n}_\infty(1)]^*$ such that $F(w) + \langle \pi, Bw \rangle \leq F(\bar{w})$ for each $w \in W$. Consider the Yosida-Hewitt decomposition $\pi = \pi^a + \pi^s$ of $\pi$ into the sum of an absolutely continuous and singular functionals $\pi^a$ and $\pi^s$ (see [27]). Let $l(\omega)$ be an integrable $\mathcal{F}_t$-measurable function representing the functional $\pi^a$, i.e. satisfying $\pi^a(y) = El(\omega)y(\omega)$ for all $y \in L^n_\infty(t)$. By repeating verbatim the arguments used in the proof of Proposition 1 in [11], one can show that the inequality $F(w) + \langle \pi, Bw \rangle \leq F(\bar{w})$ implies

$$F(u,v) + El[u - E_0 u] \leq 0 \text{ for all } (u,v) \in W_t. \tag{11}$$

Define $g := E_{t-1}l - l$. Inequality (7) follows from (11) in view of the identity $El[u - E_{t-1}u] = E[l - E_{t-1}l]u$. □

By using Theorem 2, we can immediately obtain assertion (i) of Theorem 1 as a consequence of a previous result about the existence of finite rapid paths defined in terms of Radner's duality.

*Proof of assertion (i) of Theorem 1.* By virtue of Theorem 3.1 in [8], for each finite $N \geq 1$ there exists a path $x_0, ..., x_N$ supported by an $R$-dual path. In view of Theorem 2, there exists a dual path supporting $x_0, ..., x_N$, and so the path $x_0, ..., x_N$ is rapid. □

## 3. Infinite rapid paths

*Proof of assertion (ii) of Theorem 1.* By virtue of (i), for each natural number $N$ there exist a rapid path $(x_0(N,\omega), ..., x_N(N,\omega))$ with $x_0(N,\omega) = x_0(\omega)$ and a dual path $(p_1(N,\omega), ..., p_{N+1}(N,\omega))$ supporting it. Observe that there exists a constant $C_t$ such that

$$|d| \leq C_t |c| \tag{12}$$

for all $(c,d) \in G_t^\times(\omega)$ and $t \geq 1$. This follows from the inequalities

$$|c| \max_i \hat{a}_{t-1}^i(\omega) \geq c\hat{a}_{t-1}(\omega) \geq d\hat{b}_t(\omega) \geq \gamma_t |d|$$

(see (1) and (**G.3**)). By using (12) and (2), we get $E_t |p_{t+1}(N,\omega)| \leq C_t |p_t(N,\omega)|$ (a.s.) and so $E|p_t(N,\omega)| \leq C_{t-1}...C_1 E|p_1(N,\omega)|$ for $N+1 \geq t \geq 2$. Furthermore, $|p_1(N,\omega)| \leq \delta^{-1}$ (a.s.) because $1 = p_1(N,\omega)x_0(\omega) \geq |p_1(N,\omega)|\delta$ (a.s.). Consequently, for each $t \geq 1$, we have

$$E|p_t(N,\omega)| \leq C^t \tag{13}$$

where $C^t := C_{t-1}...C_1 \delta^{-1}$ for $t \geq 2$ and $C^1 = \delta^{-1}$. Also, note that

$$|x_t(N,\omega)| \leq M^t \text{ (a.s.)}, \tag{14}$$

where $M^t := M_t...M_1 nD$, by virtue of condition (**G.2**) and the inequality $x_0(\omega) \leq De$.

We will construct by induction for each $t = 1, 2, ...$ a triplet of vector functions

$$(x_t, p_t, q_t) \in \mathcal{X}_t \times \mathcal{P}_t \times \mathcal{P}_t \tag{15}$$

and a sequence of integer-valued functions

$$t < N_1^t(\omega) < N_2^t(\omega) < ... \tag{16}$$

such that $N_k^t(\omega)$ is $\mathcal{F}_t$-measurable for all $k \geq 1$,

$$(p_t(\omega), q_t(\omega)) \in G_t^\times(\omega) \text{ (a.s.)}, \tag{17}$$

$$(x_{t-1}(\omega), x_t(\omega)) \in G_t(\omega) \text{ (a.s.)}, \tag{18}$$

$$p_t(\omega) x_{t-1}(\omega) = 1 \text{ (a.s.)}, \tag{19}$$

$$q_{t-1}(\omega) \geq E_{t-1} p_t(\omega) \text{ (a.s.)}, \tag{20}$$

$$x_t(N_k^t(\omega), \omega) \to x_t(\omega), \; p_t(N_k^t(\omega), \omega) \to p_t(\omega), \; q_t(N_k^t(\omega), \omega) \to q_t(\omega) \text{ (a.s.)}, \tag{21}$$

as $k \to \infty$, and $N_1^t(\omega), N_2^t(\omega), ...$ is a subsequence of $N_1^{t-1}(\omega), N_2^{t-1}(\omega), ...$ which can be represented as

$$N_m^t(\omega) = N_{k(m,\omega)}^{t-1}(\omega) \tag{22}$$

for some $\mathcal{F}_t$-measurable functions $0 < k(1,\omega) < k(2,\omega) < ...$ with integer values. For $t = 1$, the function $q_{t-1}(\omega) = q_0(\omega)$ involved in (20) is defined by $q_0(\omega) := \delta^{-1} e$ and the sequence $N_k^{t-1}(\omega) = N_k^0(\omega)$ involved in (22) is given by $N_k^0(\omega) := k$ for all $k \geq 1$.

Define

$$q_t(N,\omega) := E_t p_{t+1}(N,\omega) \; (0 \leq t \leq N). \tag{23}$$

Observe that for each $t = 1, 2, ...$ and $N \geq t$, we have

$$|x_t(N, \omega)| \leq M^t \text{ (a.s.)}, \ E|p_t(N, \omega)| \leq C^t, \ E|q_t(N, \omega)| \leq C^{t+1}, \tag{24}$$

$$(p_t(N, \omega), q_t(N, \omega)) \in G_t^\times(\omega) \text{ (a.s.)}, \tag{25}$$

$$(x_{t-1}(N, \omega), x_t(N, \omega)) \in G_t(\omega) \text{ (a.s.)}, \tag{26}$$

$$p_t(N, \omega)x_{t-1}(N, \omega) = 1 \text{ (a.s.)}, \tag{27}$$

$$q_{t-1}(N, \omega) = E_{t-1}p_t(N, \omega) \text{ (a.s.)}, \tag{28}$$

where the inequalities in (24) follow from (14), (13) and (23), relations (25) – (27) hold because $(p_1(N, \omega), ..., p_{N+1}(N, \omega))$ is a dual path supporting $(x_0(N, \omega), ..., x_N(N, \omega))$, and (28) is valid by the definition of $q_{t-1}(N, \omega)$.

Let us construct a triplet (15) and a sequence (16) satisfying (17) – (22) for $t = 1$. We will apply Proposition A.1 (see the Appendix) to the sequence of $3n$-dimensional $\mathcal{F}_1$-measurable random vectors

$$w_N(\omega) = (x_1(N, \omega), p_1(N, \omega), q_1(N, \omega)) \ (N \geq 1).$$

For each $N = 1, 2, ...$, these vectors satisfy condition (24) with $t = 1$, which implies that $\liminf E|w_N(\omega)| < \infty$. By virtue of Proposition A.1, there exists an $\mathcal{F}_1$-measurable vector function $w(\omega) = (x_1(\omega), p_1(\omega), q_1(\omega))$ and a sequence of $\mathcal{F}_1$-measurable integer-valued functions $1 < N_1^1(\omega) < N_2^1(\omega) < ...$ satisfying (15) and (21) with $t = 1$. Since the sets $G_t^\times(\omega)$ and $G_t(\omega)$ are closed, the relations

$$(p_1(N_k^1(\omega), \omega), q_1(N_k^1(\omega), \omega)) \in G_1^\times(\omega), \ k = 1, 2, ... \text{ (a.s.)}, \tag{29}$$

$$(x_0(\omega), x_1(N_k^1(\omega), \omega)) \in G_1(\omega), \ k = 1, 2, ... \text{ (a.s.)}, \tag{30}$$

$$p_1(N_k^1(\omega), \omega)x_0(\omega) = 1, \ k = 1, 2, ... \text{ (a.s.)} \tag{31}$$

yield in the limit (17) – (19) for $t = 1$. The equality $p_1(\omega)x_0(\omega) = 1$ (a.s.) implies $|p_1(\omega)| \leq \delta^{-1}$ (a.s.), and so $E_0 p_1(\omega) \leq \delta^{-1} e = q_0(\omega)$ (a.s.), which gives (20).

Suppose a triplet (15) and a sequence (16) satisfying (17) – (22) are constructed for some $t \geq 1$; let us construct such a triplet and a sequence for $t + 1$. Since $N_k^t(\omega) \geq t + 1$, relations (25) – (28) (with $t + 1$ in place of $t$) imply

$$(p_{t+1}(N_k^t(\omega), \omega), q_{t+1}(N_k^t(\omega), \omega)) \in G_{t+1}^\times(\omega) \text{ (a.s.)}, \tag{32}$$

$$(x_t(N_k^t(\omega), \omega), x_{t+1}(N_k^t(\omega), \omega)) \in G_{t+1}(\omega) \text{ (a.s.)}, \tag{33}$$

$$p_{t+1}(N_k^t(\omega), \omega)x_t(N_k^t(\omega), \omega) = 1 \text{ (a.s.)}, \tag{34}$$

$$q_t(N_k^t(\omega), \omega) = E_t p_{t+1}(N_k^t(\omega), \omega) \text{ (a.s.)}. \tag{35}$$

The last relation follows from the equality

$$\mathbf{1}_{\{N_k^t(\omega)=m\}} q_t(N_k^t(\omega),\omega) = \mathbf{1}_{\{N_k^t(\omega)=m\}} E_t p_{t+1}(N_k^t(\omega),\omega) \text{ (a.s.)},$$

holding for each $m = t+1, t+2, ...$ because the integer-valued random variable $N_k^t(\omega)$ is $\mathcal{F}_t$-measurable. According to the construction of the triplet $x_t(\omega), p_t(\omega), q_t(\omega)$ and the sequence $N_k^t(\omega)$, we have

$$q_t(N_k^t(\omega),\omega)) \to q_t(\omega) \text{ (a.s.)}, \ x_t(N_k^t(\omega),\omega) \to x_t(\omega) \text{ (a.s.)}. \tag{36}$$

By virtue of (35) and (36), we have

$$E_t p_{t+1}(N_k^t(\omega),\omega) \to q_t(\omega) \text{ (a.s.)}.$$

We apply the conditional version of the multidimensional Fatou lemma (see the Appendix, Proposition A.2) to the sequence of $\mathcal{F}_{t+1}$-measurable random vectors $z_k(\omega) := p_{t+1}(N_k^t(\omega),\omega)$. By virtue of this lemma, there exists a sequence $1 < k(1,\omega) < k(2,\omega) < ...$ of $\mathcal{F}_{t+1}$-measurable integer-valued functions and an $\mathcal{F}_{t+1}$-measurable vector function $p_{t+1}(\omega) \geq 0$ such that

$$z_{k(l,\omega)}(\omega) \to p_{t+1}(\omega) \text{ (a.s.)}, \tag{37}$$

$$q_t(\omega) \geq E_t p_{t+1}(\omega) \text{ (a.s.)}. \tag{38}$$

Since $q_t \in \mathcal{P}_t$, inequality (38) implies that $p_{t+1} \in \mathcal{P}_{t+1}$. By setting $n_l^{t+1}(\omega) := N_{k(l,\omega)}^t(\omega)$, we obtain a sequence $t+1 < n_1^{t+1}(\omega) < n_2^{t+1}(\omega) < ...$ of $\mathcal{F}_{t+1}$-measurable integer-valued functions such that

$$p_{t+1}(n_l^{t+1}(\omega),\omega) \to p_{t+1}(\omega), \ x_t(n_l^{t+1}(\omega),\omega) \to x_t(\omega) \text{ (a.s.)} \tag{39}$$

by virtue of (38) and (36). In view of (32) – (34), we get

$$(p_{t+1}(n_l^{t+1}(\omega),\omega), q_{t+1}(n_l^{t+1}(\omega),\omega)) \in G_{t+1}^\times(\omega) \text{ (a.s.)}, \tag{40}$$

$$(x_t(n_l^{t+1}(\omega),\omega), x_{t+1}(n_l^{t+1}(\omega),\omega)) \in G_{t+1}(\omega) \text{ (a.s.)}, \tag{41}$$

$$p_{t+1}(n_l^{t+1}(\omega),\omega) x_t(n_l^{t+1}(\omega),\omega) = 1 \text{ (a.s.)}. \tag{42}$$

Since the sequence $p(l,\omega) := p_{t+1}(n_l^{t+1}(\omega),\omega)$ converges for almost all $\omega$, it is bounded for almost all $\omega$. By virtue of (40) and (12), the sequence $q(l,\omega) := q_{t+1}(n_l^{t+1}(\omega),\omega)$ is a.s. bounded too. For $x(l,\omega) := x_{t+1}(n_l^{t+1}(\omega),\omega)$, we have $|x(l,\omega)| \leq M^{t+1}$ (a.s.) according to (14). Therefore we can apply Proposition A.1 to the sequence of $\mathcal{F}_{t+1}$-measurable random vectors $v_l(\omega) := (x(l,\omega), q(l,\omega))$. By virtue of this proposition, there exists a sequence $1 < l(1,\omega) < l(2,\omega) < ...$ of $\mathcal{F}_{t+1}$-measurable integer-valued random variables and $\mathcal{F}_{t+1}$-measurable random vectors $x_{t+1}(\omega) \geq 0, q_{t+1}(\omega) \geq 0$ for which

$$x(l(m,\omega),\omega) \to x_{t+1}(\omega), \ q(l(m,\omega),\omega) \to q_{t+1}(\omega) \text{ (a.s.)}. \tag{43}$$

Since the sets $G_{t+1}^\times(\omega)$ and $G_{t+1}(\omega)$ are closed, it follows from (40) – (43) and (39) that

$$(p_{t+1}(\omega), q_{t+1}(\omega)) \in G_{t+1}^\times(\omega) \text{ (a.s.)}, \tag{44}$$

$$(x_t(\omega), x_{t+1}(\omega)) \in G_{t+1}(\omega) \text{ (a.s.)}, \tag{45}$$

$$p_{t+1}(\omega)x_t(\omega) = 1 \text{ (a.s.)}. \tag{46}$$

We have $|x_{t+1}(\omega)| \leq M^{t+1}$ (a.s.) because $|x(l,\omega)| \leq M^{t+1}$ (a.s.). The inclusion in (44), combined with (12), yields $|q_{t+1}(\omega)| \leq C_t |p_{t+1}(\omega)|$, and so $q_{t+1} \in \mathcal{P}_{t+1}$ as long as $p_{t+1} \in \mathcal{P}_{t+1}$, which was shown above. Thus the triplet $(x_{t+1}, p_{t+1}, q_{t+1}) \in \mathcal{X}_{t+1} \times \mathcal{P}_{t+1} \times \mathcal{P}_{t+1}$ satisfies all the conditions listed in (17) – (20) (with $t+1$ in place of $t$). It remains to define the sequence $N_m^{t+1}(\omega)$, $m = 1, 2, ...$, of $\mathcal{F}_{t+1}$-measurable random integers by

$$N_m^{t+1}(\omega) = n_{l(m,\omega)}^{t+1}(\omega), \text{ where } n_l^{t+1}(\omega) = N_{k(l,\omega)}^t(\omega). \tag{47}$$

The sequence $N_m^{t+1}(\omega)$ is strictly increasing in $m$ because the sequences $N_k^t(\omega)$, $k(l,\omega)$, $n_l^{t+1}(\omega)$ and $l(m,\omega)$ are strictly increasing. We have $N_1^{t+1}(\omega) > t+1$ since $N_1^t(\omega) > t$ and $k(l,\omega) > 1$. By using formulas (43) and (39), we obtain $x_{t+1}(N_m^{t+1}(\omega), \omega) \to x_{t+1}(\omega)$, $p_{t+1}(N_m^{t+1}(\omega), \omega) \to p_{t+1}(\omega)$ and $q_{t+1}(N_m^{t+1}(\omega), \omega) \to q_{t+1}(\omega)$ (a.s.). Thus the sequence $N_m^{t+1}(\omega)$ possesses all the properties required for $t+1$ in (21) and (22) (the latter follows from (47)).

We have constructed for each $t = 1, 2, ...$, non-negative measurable vector functions $x_t(\omega)$, $p_t(\omega)$, $q_t(\omega)$ satisfying (15) and (17) – (20). Consider the sequences $(x_0(\omega), x_1(\omega), ...)$ and $(p_1(\omega), p_2(\omega), ...)$. It follows from (18) that the former is a path. The latter is a dual path by virtue of relation (17), the inequality $q_t(\omega) \geq E_t p_{t+1}(\omega)$ (a.s.) and the fact that the cone $G_t^\times(\omega)$ contains with any pair of vectors $(c, d)$ any pair $(c, d')$ with $0 \leq d' \leq d$. Finally, $(p_1(\omega), p_2(\omega), ...)$ supports $(x_0(\omega), x_1(\omega), ...)$ by virtue of (19). □

## Acknowledgment

Financial support from the grant NSF DMS-0505435, the State of Missouri Research Board, the University of Missouri-Columbia Research Council, and the Manchester School Visiting Fellowship Fund is gratefully acknowledged.

**Appendix A:**

Let $(\Omega, \mathcal{F}, P)$ be a probability space and $w_N(\omega)$, $\omega \in \Omega$, $N = 1, 2, ...$, a sequence of random vectors in $\mathbb{R}^n$.

**Proposition A.1.** *If* $\liminf |w_N(\omega)| < \infty$ *a.s. (which is so, in particular, when* $\liminf E|w_N(\omega)| < \infty$*), then there exists a sequence of integer-valued random variables* $N_1(\omega) < N_2(\omega) < ...$ *and a random vector* $w(\omega)$ *such that*

$$\lim w_{N_k(\omega)}(\omega) = w(\omega) \ (a.s.)$$

*and*

$$E|w(\omega)| \leq \liminf E|w_N(\omega)|. \tag{A1}$$

*Proof.* Define $\xi(\omega) = \liminf |w_N(\omega)|$. Since $\xi(\omega) < \infty$ (a.s.), we can select a sequence of natural numbers $N_1(\omega) < N_2(\omega) < ...$ such that $\lim |w_{N_k(\omega)}(\omega)| = \xi(\omega)$ (a.s.) and $w_{N_k(\omega)}(\omega) \to w(\omega)$ (a.s.), where $w(\omega)$ is a vector with $|w(\omega)| = \xi(\omega)$. A measurable selection argument based on Aumann's theorem (see e.g. [1], Appendix I, Corollary 3) makes it possible to choose $N_k(\omega)$, $k = 1, 2, ...$, and $w(\omega)$ in a measurable way. Inequality (A1) follows from the Fatou lemma and the equality $|w(\omega)| = \liminf |w_N(\omega)|$. □

**Proposition A.2.** *Let* $\mathcal{G}$ *be a sub-$\sigma$-algebra of* $\mathcal{F}$. *If* $w_N(\omega)$, $\omega \in \Omega$, $N = 1, 2, ...$, *are random vectors with values in* $\mathbb{R}_+^n$ *and the conditional expectations* $E[w_N(\omega)|\mathcal{G}]$ *are finite and converge a.s. to a random vector* $z(\omega)$, *then exists a sequence of integer-valued random variables* $1 < N_1(\omega) < N_2(\omega) < ...$ *and a random vector* $w(\omega)$ *such that* $\lim w_{N_k(\omega)}(\omega) = w(\omega)$ *(a.s.) and*

$$E[w(\omega)|\mathcal{G}] \leq z(\omega) \ (a.s.).$$

This result is a version of the multidimensional Fatou lemma for conditional expectations. For a proof of its conventional version dealing with unconditional expectations (corresponding to the case when the $\sigma$-algebra $\mathcal{G}$ is trivial) see Schmeidler [24]. A lemma closely related to Schmeidler's one was obtained by Olech [17]. For applications of such results to mathematical economics and optimal control theory, see Hildenbrand [15] and Olech [18].

*Proof of Proposition A.2.* Let us consider the sequence of random vectors $w^\infty(\omega) := (w_1(\omega), w_2(\omega), ...)$ as a random element of the standard measurable space $(X^\infty, \mathcal{B}^\infty) := (X, \mathcal{B}) \times (X, \mathcal{B}) \times ...$, where $X := \mathbb{R}_+^n$ and $\mathcal{B}$ is the Borel $\sigma$-algebra on $X$. Let $\pi(\omega, dx^\infty)$ be the regular conditional distribution of $w^\infty(\omega)$ given the $\sigma$-algebra $\mathcal{G}$ (see Shiryaev [25], Theorem II.7.5). Denote by $x_N(x^\infty)$ the $N$th element of the sequence $x^\infty = (x_1, x_2, ...)$. We have

$$E[w_N(\omega)|\mathcal{G}] = \int \pi(\omega, dx^\infty) x_N(x^\infty) \to z(\omega)$$

for all $\omega$ belonging to some $\mathcal{G}$-measurable set $\Omega_1 \subseteq \Omega$ of measure one. By virtue of the conventional (unconditional) version of the Fatou lemma in several dimensions, for each $\omega \in \Omega_1$ there exists a $\mathcal{B}^\infty$-measurable vector function $v^\omega(x^\infty)$ such that $\int \pi(\omega, dx^\infty) v^\omega(x^\infty) \leq z(\omega)$ and $v^\omega(x^\infty) \in \mathrm{Ls}\,(x^\infty)$ for $\pi(\omega, \cdot)$-almost all $x^\infty$, where $\mathrm{Ls}\,(x^\infty)$ is the set of limit points of the sequence $x^\infty = (x_1(x^\infty), x_2(x^\infty), ...)$.

We will use the following fact. There exists a function $\psi : [0,1] \times X^\infty \to \mathbb{R}^n_+$ jointly measurable with respect to $\mathcal{B}[0,1] \times \mathcal{B}^\infty$ (where $\mathcal{B}[0,1]$ is the Borel $\sigma$-algebra on $[0,1]$) and possessing the following property. For each finite measure $\mu$ on $\mathcal{B}^\infty$ and each $\mathcal{B}^\infty$-measurable function $f : X^\infty \to \mathbb{R}^n_+$, there exists $r \in [0,1]$ such that $\psi(r, x^\infty) = f(x^\infty)$ for $\mu$-almost all $x^\infty \in X^\infty$. This result establishes the existence of a "universal" jointly measurable function parametrizing all equivalence classes of measurable functions $X^\infty \to \mathbb{R}^n_+$ with respect to all finite measures: any such class contains a representative of the form $\psi(r, \cdot)$, where $r$ is some number in $[0,1]$. The result (extending Natanson [16], Chapter 15, Section 3, Theorem 4) follows from Theorem AI.3 in [12] and the fact that all uncountable standard Borel spaces are isomorphic to the segment $[0,1]$ equipped with the Borel $\sigma$-algebra.

For each $\omega \in \Omega$, consider the set $U(\omega)$ of those $r \in [0,1]$ for which the function $\psi(r, \cdot)$ satisfies

$$\int \pi(\omega, dx^\infty) \psi(r, x^\infty) \leq z(\omega), \; \psi(r, x^\infty) \in \mathrm{Ls}\,(x^\infty) \text{ for } \pi(\omega, \cdot)\text{-almost all } x^\infty. \tag{A2}$$

Observe that for $\omega \in \Omega_1$ the set $U(\omega)$ is not empty because it contains that element $r \in [0,1]$ for which $\psi(r, x^\infty) = v^\omega(x^\infty)$ $\pi(\omega, \cdot)$-almost everywhere. Further, the set of pairs $(\omega, r)$ satisfying (A2) is $\mathcal{G} \times \mathcal{B}[0,1]$-measurable because $\pi(\omega, dx^\infty)$ is a conditional distribution given $\mathcal{G}$, the function $\psi(r, x^\infty)$ is $\mathcal{B}[0,1] \times \mathcal{B}^\infty$-measurable and the second constraint in (A2) can be written as

$$\int \pi(\omega, dx^\infty) F(r, x^\infty) = 1, \tag{A3}$$

where $F(r, x^\infty)$ is the indicator function of the set

$$\{(r, x^\infty) : \psi(r, x^\infty) \in \mathrm{Ls}\,(x^\infty)\} \in \mathcal{B}[0,1] \times \mathcal{B}^\infty.$$

The last inclusion follows from the fact that $y \in \mathrm{Ls}\,(x^\infty)$ if and only if for each $m = 1, 2, ...$ and $k = 1, 2, ...$ there exists $K \geq k$ such that $|y - x_K(x^\infty)| < 1/m$.

By virtue of Aumann's measurable selection theorem, there exists a $\mathcal{G}$-measurable function $r(\omega)$ such that $r(\omega) \in U(\omega)$ (a.s.). Define $w(\omega) := \psi(r(\omega), w^\infty(\omega))$. Since $\pi(\omega, dx^\infty)$ is the conditional distribution of $w^\infty(\omega)$ given $\mathcal{G}$, we have

$$E[w(\omega)|\mathcal{G}] = E[\psi(r(\omega), w^\infty(\omega))|\mathcal{G}] = \int \pi(\omega, dx^\infty) \psi(r(\omega), x^\infty) \leq z(\omega) \text{ (a.s.)}.$$

Furthermore, $w(\omega) \in \mathrm{Ls}\,(w^\infty(\omega))$ (a.s.) because this inclusion is equivalent to the equality $F(r(\omega), w^\infty(\omega)) = 1$ (a.s.) and

$$EF(r(\omega), w^\infty(\omega)) = E\{E[F(r(\omega), w^\infty(\omega))|\mathcal{G}]\} =$$

$$E \int \pi(\omega, dx^\infty) F(r(\omega), x^\infty) = 1$$

by virtue of (A2) and (A3). Since $w(\omega) \in \mathrm{Ls}\,(w^\infty(\omega))$ (a.s.), it follows from Aumann's measurable selection theorem that there exist integer-valued random variables $1 < N_1(\omega) < N_2(\omega) < ...$ such that $\lim w_{N_k(\omega)}(\omega) = w(\omega)$ (a.s.). □

**Proposition A.3.** *Let $f(\omega, x)$ be a real-valued jointly measurable function of $\omega \in \Omega$, $x \in \mathbb{R}^n$ and $G(\omega) \subseteq \mathbb{R}^n$ a non-empty closed set measurably depending on $\omega$ such that for $x \in G(\omega)$ we have $|f(\omega, x)| \leq$* $K(\omega)$, where $EK(\omega) < \infty$. Let $\bar{\xi}(\omega)$ be a measurable vector function such that $\bar{\xi}(\omega) \in G(\omega)$ (a.s.). Then the following assertions are equivalent. (i) The inequality $Ef(\omega, \xi(\omega)) \leq Ef(\omega, \bar{\xi}(\omega))$ holds for all measurable functions $\xi(\omega)$ satisfying $\xi(\omega) \in G(\omega)$ (a.s.). (ii) With probability one,

$$f(\omega, \bar{\xi}(\omega)) = \max_{x \in G(\omega)} f(\omega, x). \tag{A4}$$

The proposition can be obtained by using Aumann's measurable selection theorem (see the reference above).

Let $D_1$ and $D_2$ be Banach spaces, $W$ a convex subset of $D_1$, $F : W \to (-\infty, +\infty)$ a concave functional and $B : D_1 \to D_2$ a continuous linear operator such that the image $H$ of the space $D_1$ under the mapping $B : D_1 \to D_2$ is a closed subspace of $D_2$. Assume that the functional $F(w)$ attains its maximum on the set $\{w \in W : Bw = 0\}$ at some point $\bar{w}$.

**Proposition A.4.** *Let the following condition hold.*

*(\*) The functional $F(\cdot)$ is bounded below on some subset $\tilde{W}$ in $W$ such that the interior of $B(\tilde{W})$ in the space $H$ contains zero.*

*Then there exists a continuous linear functional $\pi$ on the space $D_2$ such that*

$$F(w) + \langle \pi, Bw \rangle \leq F(\bar{w}) \tag{A5}$$

*for each $w \in W$.*

This result is a version of the Kuhn-Tucker theorem for concave optimization problems with linear equality constraints. Its proof is similar to that of Proposition 1 in [9].